\documentclass[runningheads]{llncs}
\usepackage[T1]{fontenc}
\usepackage[utf8]{inputenc}
\usepackage[datamodel=software]{biblatex}
\usepackage{software-biblatex}
\usepackage{amsmath,amssymb,amsfonts}
\usepackage{hyperref}
\usepackage{dcolumn}
\newcolumntype{d}{D{.}{.}{-1}}

\newcommand{\point}[1]{(#1)}

\renewcommand{\footnotemark}{}

\renewcommand {\geq}{\geqslant}
\renewcommand {\leq}{\leqslant}
\newcommand {\parigp}{\textsc {Pari/Gp}}
\newcommand {\cm}{\textsc {Cm}}
\newcommand{\legendre}[2]{\genfrac {(}{)}{1pt}{}{#1}{#2}}
\newcommand {\Z}{\mathbb {Z}}
\newcommand {\Q}{\mathcal {Q}}
\newcommand {\Ot}{\tilde {O}}
\newcommand {\Omegat}{\tilde {\Omega}}
\newcommand {\Thetat}{\tilde {\Theta}}
\newcommand {\Dmax}{D_{\mathrm {max}}}
\newcommand {\hmax}{h_{\mathrm {max}}}
\newcommand {\pmax}{p_{\mathrm {max}}}

\title{FastECPP over MPI}
\author{Andreas Enge}
\institute {INRIA, Université de Bordeaux, CNRS, CANARI,
33400 Talence, France \newline
\url{https://enge.math.u-bordeaux.fr/} \newline
\url{andreas.enge@inria.fr}}
\date{May 1st, 2024}

\begin{document}
\maketitle

\begin {abstract}
The FastECPP algorithm is currently the fastest approach to prove the
primality of general numbers, and has the additional benefit of creating
certificates that can be checked independently and with a lower complexity.
This article shows how by parallelising over a linear number of cores,
its quartic time complexity becomes a cubic wall-clock time complexity;
and it presents the algorithmic choices of the FastECPP implementation in
the author's \cm\ software \\
\centerline {\url {https://www.multiprecision.org/cm/}}
which has been written with massive parallelisation over MPI in mind,
and which has been used to establish a new primality record for the
``repunit'' $(10^{86453} - 1) / 9$.
\end {abstract}

\section {FastECPP and its complexity}

Since its inception, ECPP has become the fastest practical algorithm for
proving the primality of arbitrary numbers~$N$; additionally, it creates
a certificate that can be verified independently and in less time.
It is based on the key observation that if one can find an elliptic curve
over $\Z / N \Z$ and a point on the curve of sufficiently large order~$N'$,
then the primality of~$N'$ implies the primality of~$N$.
One step of the certificate creation process consists of searching for a
suitable value of $N'$, smaller than $N$ by at least one bit, the elliptic
curve and the point; then it continues recursively with~$N'$.
If $L = \lceil \log_2 (N) \rceil$, this requires $O (L)$ steps.
Verification goes through these $O (L)$ steps again, but checking their
correctness is much faster than finding them.

To find a suitable curve, Goldwasser and Kilian in their original
algorithm~\cite {GoKi86} use point counting on random curves, which
results in a (heuristic, under the assumption that all occurring numbers
behave randomly with respect to the sizes of their prime factors) time
complexity of $\Ot (L^6)$ using asymptotically fast arithmetic, where the
$\Ot$ notation hides logarithmic factors.
Atkin and Morain observe in~\cite {AtMo93} that one can find suitable
elliptic curves using the complex multiplication method, which lowers
the heuristic complexity to $\Ot (L^5)$; this is now known as the ECPP
algorithm.

The FastECPP version, attributed to Shallit in \cite[\S5.10]{LeLe90} and
worked out in~\cite{FrKlMoWi04,Morain07}, improves a bottleneck, lowering
the heuristic complexity to $\Ot (L^4)$.

In this section, we present the FastECPP algorithm and its heuristic
analysis following~\cite{Morain07}; for a comprehensive overview of
primality testing and proving algorithms, see~\cite {Morain04}.
Basic ECPP can essentially be recovered from FastECPP by dropping
substep~\point {1} of the first phase below and by computing the square
roots of the discriminants in substep~\point {2} one by one.
Since this is slower in all cases and does not even substantially simplify
the implementation, the basic variant is only of historical interest.

\paragraph {First phase.}

One step of the certificate creation process determines a probable prime
number $N'$ that is smaller than $N$ and then recursively creates a
certificate for $N'$.
Each of the steps consists of the following substeps.
\begin {enumerate}
\item
Write down the set $\Q = \{ q_1^*, q_2^*, \ldots \}$
of $\Thetat (L)$ smallest ``signed primes'' with
$\legendre {q^*}{N} = 1$, where
$q^* = q$ if $q$ is a prime that is $1 \bmod 4$,
$q^* = -q$ if $q$ is a prime that is $3 \bmod 4$,
or $q^* \in \{ -4, \pm 8 \}.$
Compute their square roots modulo $N$, in time $\Ot (L^2)$
each or $\Ot (L^3)$ altogether.
\item
Fundamental quadratic discriminants are exactly the products of signed
primes (without multiplicities); create a set of $\Thetat (L^2)$ negative
discriminants~$D$ with $|D| \leq \Dmax \in \Ot (L^2)$
and their square roots modulo~$N$ as products of the
precomputed roots (products of two elements of~$\Q$ are enough).
Try to solve the Pell equation $4 N = t^2 - v^2 D$ for all~$D$ using
Cornachia's algorithm, in time $\Ot (L)$ per problem or
$\Ot (L^3)$ altogether. The success probability is of the order
of $1 / \sqrt {|D|} \in \Omegat (1 / L)$, so $\Thetat (L)$ values $t$ are
expected to be obtained; then if $N$ is prime, there are elliptic curves
modulo~$N$ with complex multiplication by~$D$ and $m = N + 1 \pm t$ points.
\item
Trial factor the $m$ as $m = c N'$, where the cofactor $c \geq 2$ is
$B$-smooth for some bound~$B$ (for instance, $B = 2$), and all primes
dividing $N'$ are larger than~$B$; each factorisation takes $\Ot (L^2)$
for $B \in \Ot (L^3)$ for a total of $\Ot (L^3)$ (more details are given
in \S\ref {sec:implementation}).
\item
Test the $N'$ for primality, in time $\Ot (L^2)$ each for a total of
$\Ot (L^3)$; the expected number of remaining primes $N'$ is
in $\Theta (1)$.
\end {enumerate}
So each step, of which there are $O (L)$, has a complexity of $\Ot (L^3)$.

\paragraph {Second phase.}

For each of the $O (L)$ steps of the first phase, carry out the
following substeps:
Construct the class polynomial of degree
$h \in \Ot \left(\sqrt {|D|} \right)$, where
$h$ is the class number of~$D$, in time
$\Ot (|D|) \subseteq \Ot (L^2)$~\cite {Enge09}.
Find a root modulo $N$ in time $\Ot (L^3)$ and write down the two
corresponding CM elliptic curves; take random points on the curves
in time $\Ot (L^2)$ for the square roots yielding the $Y$-coordinates and
multiply them by the cofactors $c$ to obtain a point of prime order~$N'$.
So this phase also has $O (L)$ independent steps of complexity
$\Ot (L^3)$ each.

\smallskip
Notice that all substeps above admit a trivial parallelism of $\Theta (L)$;
so using $\Omegat (L)$ cores, the wall-clock time complexity of FastECPP
becomes $\Ot (L^3)$.
Some of the substeps could be further parallelised to $\Theta (L^2)$ cores,
but it is probably not possible to improve the complexity of a modular
square root or a primality test beneath $\Omegat (L^2)$ even in parallel.

The final certificate is of bit size $O (L^2)$ and can be verified in
sequential time $\Ot (L^3)$, or in wall-clock time $\Ot (L^2)$
on $\Omegat (L)$ cores.

\section {Implementation choices}
\label {sec:implementation}

The only previously available free implementation of the FastECPP algorithm
was written by Jared Asuncion within \parigp~\cite {parigp-2.15.4} as the
\texttt {primecert} function, with numbers of around 1000 digits in mind.
The \cm\ software \cite {cm-0.4.3}, on the other hand, encapsulated over two
decades of my research on algorithms for CM of elliptic curves; so I decided
to add first a sequential FastECPP implementation and eventually an MPI
version tuned for massive parallelism and with record computations in mind.
Both are available in \cm\ since version 0.4.0 as the binaries
\texttt {ecpp} and \texttt {ecpp-mpi}, respectively.
This section documents the choices made in the implementation and the
parameter selection depending on the input size~$L$ and the available
computing power.
Indeed the MPI version dedicates the main process to the coordination of
the computation and uses the $w$ remaining ones to do the actual work.
It adapts to the number of available cores by choosing parameters
depending on~$w$; so also the final certificate depends on~$w$, see~\S5.

\paragraph {First phase. Substep \point {0}.}

Inspired by \parigp, the class numbers of all discriminants up to $\Dmax$
are precomputed in time $\Ot (\Dmax^{1.5})$ by looping over all quadratic
forms of discriminant in the desired range.
Also, a number of prime products used for trial division are precomputed.
Both are parallelised over all $w$ workers, and the results can be stored
on disk for reuse over several runs.

The following substeps \point {1} to \point {4} may need to be repeated
over several rounds in the unlucky case that no candidate $N'$ remains
at the end.

\paragraph {Substep \point {1}.}
Let $\Dmax = \min \left( 2^{35}, \max (2^{20}, L^2/2) \right)$
(or, more precisely, the smallest power of~$2$ not below this),
$\hmax = 100000$ and
$\pmax \! = \max (29, \lfloor L/2^{10} \rfloor)$.
The values of $\hmax$ and $\pmax$ are intended to speed up the second
phase:
Only discriminants with class number at most $\hmax$ are considered,
which is an upper bound on the degree of class polynomials constructable
in reasonable time;
and the maximal prime dividing the class number is bounded by $\pmax$,
so that after expressing the class field as a tower of prime degree
extensions using \cite {EnMo03}, it is enough to determine roots of
polynomials of degree at most $\pmax$.

The main process computes the set $\Q$ of $k w$ smallest signed primes
and the discriminants divisible by up to $7$ primes from~$\Q$ under the
additional restrictions above.
The integer $k$ is chosen minimally such that the expected number of
remaining~$N'$ after substep~\point {4}, obtained using the formula for $s$
in \cite [\S4]{FrKlMoWi04}, the formula in the last line of
\cite [\S2]{Morain07} for the success probability of solving the
Pell equation and the precomputed class numbers, is at least~$3$
in the first round (to allow for some choice), or~$1$ in later rounds (to
avoid computing too many square roots).
Then each worker computes $k$ square roots modulo~$N$.

\paragraph {Substep \point {2}.}

The implementation evenly distributes \textit {all} discriminants among
the workers to solve the Pell equation. If the expected number of $N'$ is
much larger than~$1$, this may lead to unnecessary work being executed,
but on the other hand, it may also lead to more choice and thus a smaller
$N'$ and ultimately fewer steps. In the absence of a clear optimality
criterion, the design decision was to favour smaller certificates over
faster running times; see the comparison with \parigp\ in
\S\ref {sec:comparison}.

\paragraph {Substep \point {3}.}

Let $P$ be the product of all primes up to~$B$;
then $|\log P / B - 1| \in o (1)$ by \cite [Th.~4]{RoSc62}.
The numbers $m_i$ of size $\Theta (L)$ obtained in the previous substep need
to be written as $m_i = c_i N_i'$ with $c_i$ having only prime factors
dividing~$P$, and $N_i'$ coprime to~$P$.
The implementation follows \cite [\S4]{FrKlMoWi04} and proceeds in batches
to first compute all the $P \bmod m_i$. In a first step, it constructs
bottom-up a binary tree with the $m_i$ on the leaves and each inner node
being the product of its two descendants; then it replaces the root
$M = m_1 m_2 \cdots$ by $P \bmod M$, and top-down each node by its parent
modulo the node, ending up with $P \bmod m_i$ on the leaves. We may
split the $m_i$ into batches, handled independently, such that $M \leq P$.
Then each batch contains $O (\log_2 P / L)$ numbers $m_i$, and the tree
computation takes time $\Ot (\log P)$. If we have $\Omegat (\log_2 P / L)$
numbers $m_i$, then the amortised cost per $m_i$ is in $\Ot (L)$.
But we may also use larger values of~$B$, such that $\log P \gg \log M$;
the sequential version \texttt {ecpp} uses $B \in \Theta (L^3)$
for $\Thetat (L)$ values~$m_i$ with an amortised cost of $\Ot (L^2)$
per value.
Now we compute the square-free part of $c_i$ as
$c_i' = \gcd (m_i, P \bmod m_i)$ in time $\Ot (L)$,
the square part of $c_i$ as $c_i'' = \gcd (m_i / c_i', c_i')$
and so on. This gcd part usually has a total complexity of $\Ot (L)$,
but in unlikely cases (for instance if $m_i = 2^L$) may require
$\Ot (L^2)$. (A complexity of $\Ot (L)$ could be achieved by using the
algorithm of \cite {Bernstein04}, which is expected to be slightly
slower on average.)

Notice that the tree-based approach requires $\Thetat (\log P)$ of memory,
and that Moore's law acts similarly on the number of cores and on the
memory; indeed one has observed over the last decades that the memory
per core has increased only very moderately.
Otherwise said, the total memory available over $\Theta (L)$ cores
is of the order of $\Ot (L)$, which implies that $B \in \Ot (L)$:
So asymptotically in the parallel setting, the effect of batch trial
division vanishes.

The MPI implementation is parallelised in two dimensions, with respect
to the~$m_i$ and to the primes in~$P$. It splits the $m_i$ into
$b \approx 16$ batches, each of which is sent to an MPI communicator with
$w' = \lfloor w/b \rfloor$ workers. Then worker~$i$ handles the product
$P_i$ of primes between $(i-1) \cdot 2^{29}$ and $i \cdot 2^{29}$, 
so that the effective smoothness bound becomes $B = w' \cdot 2^{29}$.
By the already cited \cite [Th.~4]{RoSc62}, each of the $P_i$ has about
the same bit size of
$\log_2 (P_i) \approx 2^{29} / \log (2) \approx 775 \cdot 10^6$,
so that each level of the tree requires
about $97$~MiB of memory. In the $86453$ digit record, there are
up to $2^{29} / (86453 \, \log (10)) \approx 2700$ values $m_i$ on the
leaves, and thus $\left\lceil \log_2 (2700) \right\rceil + 1 = 13$
levels, so the total memory requirements are about $1.3$~GiB per core.

\paragraph {Substep \point {4}.}

The implementation performs Miller--Rabin tests in batches of~$w$ numbers,
starting with the smallest ones. So of all the non-smooth parts computed
in substep~\point {3}, the smallest prime one will be retained as $N'$.

\section {Details of the 86453 digit record}

The second phase of FastECPP is built upon a number of research results
that were already implemented in the CM software~\cite {cm-0.4.3}.
It computes class polynomials by complex approximations as described
in~\cite {Enge09}. Optimal class invariants are chosen derived from Weber
functions~\cite {Schertz76},
simple~\cite {EnMo14} or double eta quotients~\cite {EnSc04},
including cases where it is enough to compute lower-degree subfields of
the class field~\cite {EnSc13}. The evaluation of modular functions,
which is the most important part of the class polynomial computation,
is optimised following~\cite {Enge09,EnHaJo18}.
To ease the step of factoring class polynomials modulo primes, the class
fields are then represented as towers of cyclic Galois extensions of
prime degree following~\cite {EnMo03}.
The software relies on a number of libraries from the GNU project, notably
GMP~\cite {gmp-6.2.1}, MPFR~\cite {mpfr-4.1.0} and MPC~\cite {mpc-1.2.1},
and on \parigp~\cite {parigp-2.15.4} for computations with class groups
and Flint~\cite {flint-2.9.0} for root finding modulo a prime.

Computations have been carried out on the PlaFRIM cluster in Bordeaux,
\url {https://www.plafrim.fr/}, to prove primality of the ``repunit''
$(10^{86453} - 1) / 9$.
The certificate is available in \parigp\ format at \\
\centerline {
\url {https://www.multiprecision.org/downloads/ecpp/cert-r86453.bz2}.
}
An independent verification of the certificate has been carried out
by the factordb server, see \\
\centerline {
\url {http://factordb.com/index.php?id=1100000000046752372}
}
The first phase computed 2980 steps in about 103 days of wall-clock
time and 383 years of CPU time, in several runs with 759 to 2639 cores
depending on machine availability.
Of this CPU time, about
13\% were devoted to the computation of square roots in
substep~\point {1},
47\% to solving Pell equations in substep~\point {2},
10\% to batch trial factoring in substep~\point {3} and
30\% to primality tests in substep~\point {4}.
The largest occurring discriminant in absolute value was $-34223767071$,
the largest prime~$q^*$ was $240869$ of the discriminant $-3329532187$,
the largest prime of a class number was 277, appearing for 16 different
discriminants.
The effective trial factor bound $B = w' \cdot 2^{29}$ varied depending
on the number $w$ of cores assigned to a run, with $43 \leq w' \leq 172$.

The second phase was run on a machine with 96 cores and 1~TB of memory
and took about 25 years of CPU time (wall-clock time was lost),
that is, only about 6\% of the total CPU time of both phases.
Inside the phase, 2.5\% of the time was devoted to computing tower
representations for class fields,
2.8\% to verifying orders of points on elliptic curves,
and close to 95\% to finding roots of class polynomials.
Depending on the largest factor of the class number, which determines
the degrees of the polynomials, running times for this dominant
step vary considerably.
The longest step took close to 42 days for factoring a class
polynomial for the discriminant $-2083578323$ of class number
$11920 = 2^4 \cdot 5 \cdot 149$ for an intermediate prime of
$82089$ digits;
this would also have been the wall-clock time for this phase
had it been run sufficiently in parallel.

Verification of the certificate took about 48 hours using \parigp\
on the same machine with 96 cores.

\section {Sizes of smooth parts}
\label {sec:size}

It is shown in \cite [\S4]{FrKlMoWi04} that the probability of obtaining
a prime quotient after removing the $B$-smooth part from an $L$-bit number
is asymptotically $e^\gamma \log_2 B / L$, where $\gamma = 0.577\ldots$
is the Euler--Mascheroni constant, and that in this case the expected
number of binary digits of the $B$-smooth part is $\log_2 B \in \Ot (1)$
for $B$ polynomial in~$L$. The following computations give more precise
estimates. They rely on approximations from analytic number theory for
which very tight error bounds are available, for instance
in~\cite {RoSc62}. As a consequence, the estimates are relevant for the
sizes of numbers we consider, say for $B \geq 2^{29}$ and $L$ of the order
of a few thousands.
To simplify the exposition, however, only the main terms are given,
using $\sim$ to denote asymptotic equality up to lower order terms.

The probability that a number between $2^{L-1}$ and $2^L$ is a $B$-smooth
number times a prime is given by
\begin {eqnarray*}
\frac {1}{2^{L-1}}
\sum_{f B\text {-smooth}} \enspace
\sum_{p \text { prime, } 2^{L-1}/f < p \leq 2^L/f}
1
& \sim &
\frac {1}{2^{L-1}}
\sum_{f B\text {-smooth}}
\frac {2^{L-1}}{f \log (2^L/f)} \\
& \sim & \frac {1}{\log (2^L)}
\sum_{f B\text {-smooth}} \frac {1}{f},
\end {eqnarray*}
using the prime number theorem and $\log f \in o (L)$.

For $0 < \alpha \leq 1$, all $f \leq B^\alpha$ are $B$-smooth, and their
contribution to the sum is
$\sum_{f \leq B^\alpha} (1/f) \sim \alpha \log B$
using the main term of the harmonic number.
So the conditional probability that $f \leq B^\alpha$
for an $L$-bit number assumed to be a $B$-smooth number~$f$ times a prime
(for $L$ tending to infinity and $B$ growing slowly) is
$\alpha / e^\gamma$,
and the probability density function for this event with respect to
$\alpha$ is the constant $1 / e^\gamma$.
In particular, $1 / e^\gamma \approx 56\%$ of the primes $N'$ remaining
after substep~\point {4} gain less than $\log_2 B$ bits compared to~$N$.

For $1 < \alpha \leq 2$, numbers $f$ with $B < f \leq B^\alpha$ that are
\textit {not} $B$-smooth are divisible by exactly one prime $P > B$,
and what remains is $B$-smooth. So their contribution to the sum is
\[
\begin {split}
& \sum_{B < f \leq B^\alpha} \frac {1}{f} -
\sum_{B < P \leq B^\alpha, \, P \text { prime}} \frac {1}{P}
   \sum_{f \leq B^\alpha / P} \frac {1}{f} \\
& \sim
(\alpha - 1) \log B -
\sum_{B < P \leq B^\alpha, \, P \text { prime}} \frac {1}{P}
   (\alpha \log B - \log P) \\
& \sim
(\alpha - 1) \log B -
\alpha \log B \big( \log \log (B^\alpha) - \log \log (B) \big) +
\big( \log (B^\alpha) - \log (B) \big) \\
& = (2 \alpha - 2 - \alpha \log \alpha) \log B
\end {split}
\]
using \cite [Thm.~5 and 6]{RoSc62}.
So the probability density function is
$(1 - \log \alpha) / e^\gamma$ for this range of $\alpha$.
In particular, $(2 - 2 \log 2) / e^\gamma \approx 34\%$ of the primes $N'$
remaining after substep~\point {4} gain between $\log_2 B$ and $2 \log_2 B$
bits compared to~$N$, and only $10\%$ of the numbers gain more.

For higher values of $\alpha$, exact computations become unwieldy,
but the previous computation still yields a \textit {lower bound}
(since non-smooth numbers with more than one large prime factor
are subtracted multiple times)
on the contribution of $B$-smooth $f \leq B^\alpha$ to the sum as
\[
\sum_{f \leq B^\alpha} \frac {1}{f} -
\sum_{B < P \leq B^\alpha, \, P \text { prime}} \frac {1}{P}
   \sum_{f \leq B^\alpha / P} \frac {1}{f}
\sim
(2 \alpha - 1 - \alpha \log \alpha) \log B.
\]
The function reaches its maximum $(e - 1) \log B$ for $\alpha = e$.
Otherwise said, the probability that an $N'$ gains more than
$2.72 \log_2 B$ bits over $N$ is less than
$1 - (e - 1) / e^\gamma \leq 3.6\%$.

\section {Comparison}
\label {sec:comparison}

As seen in the previous section, the probability that a candidate for $N'$
gains much more than $\log_2 B$ bits is rather small, so maybe spending
almost half of the time for solving Pell equations in the record to obtain
more prime candidates for~$N'$ is not optimal. On the other hand,
probabilities of gaining more digits in one step (and thus requiring fewer
steps and less time) are amplified by the maximum statistics. For instance,
the expected number of prime candidates for $N'$ in the first step of the
record was~$8.9$; so the probability of gaining at least $2 \log_2 B$ bits
was about $1 - \big( (3 - 2 \log 2) / e^\gamma \big)^{8.9} \approx 58\%$,
and in fact $77 = 2.2 \cdot \log_2 B$ bits were gained.
Precisely, with $w = 1135$ workers, $41440753$ discriminants were
considered in substep~\point {2}, $40788$ curve cardinalities were
trial factored in substep~\point {3}, and $3405$ primality tests were
carried out in substep~\point {4} before finding a prime~$N'$ (which is
consistent with the expected number of $8.9$ primes in a sample of
$40788$).

The following table compares the wall-clock times (in minutes) and the
lengths of the certificates between \parigp~\cite {parigp-2.15.4} (which
is parallelised using threads) and \cm~\cite {cm-0.4.3} for the first
prime after $10^n$ for different values of~$n$.
The first two column blocks show figures for the parallel versions on
a machine with 128 cores. The last column block provides the results
for the serial code in \cm\ on the same type of machine.
The initial smoothness bound $B$ decreases (down to $2^{20}$) for \parigp\
as the $N'$ become smaller, and remains fixed for \cm\ over the course
of the algorithm.
Notice that the storage requirement and the verification time of
a certificate for a number of given size are proportional to the
number of steps, so that shorter certificates are preferable.

\begin {center}
\begin {tabular}{r|rrd|rrd|rrd}
\multicolumn {1}{c|}{$n$}
      & \multicolumn {3}{c|}{\parigp}
      & \multicolumn {3}{c|}{\cm\; \texttt {ecpp-mpi}}
      & \multicolumn {3}{c}{\cm\; \texttt {ecpp}} \\
      & $\log_2 B$ & \#steps & \multicolumn {1}{c|}{time}
      & $\log_2 B$ & \#steps & \multicolumn {1}{c|}{time}
      & $\log_2 B$ & \#steps & \multicolumn {1}{c}{time} \\
\hline
 1000 &  24 & 131    &   0.30 &  32 &   33 &   5.2 & 22 &  88 &   0.80 \\
 2000 &  26 & 220    &   3.6  &  32 &   76 &  15   & 25 & 157 &  15    \\
 4000 &  30 & 344    &  41    &  32 &  166 &  47   & 28 & 274 & 210 \\
 5000 &  30 & 399    &  92    &  32 &  204 &  51   & 29 & 342 & 510 \\
10000 &  30 & 740    & 820    &  32 &  444 & 220   & 32 &     &
\end {tabular}
\end {center}

Let us first compare the parallel implementations.
The \cm\ certificates are considerably shorter with an advantage that
decreases as \parigp\ chooses larger smoothness bounds~$B$.
The number of steps increases roughly linearly as expected
(slightly more for \cm, slightly less for \parigp).
One notices that the average gain of bits per step,
computed as $\log_2 (10^n)$ divided by the number of steps,
is above $2.3 \, \log_2 B$ for \cm, illustrating the effect of the
maximum statistic.

Somewhat surprisingly, neither of the two implementations shows the
quartic asymptotic running time for a fixed number of cores, which may
be due to ``overparallelisation'' and thus less than optimal CPU times
for the smaller instances: \cm\ spends $78\%$ of the wall-clock
time on trial factorisation for 5000 digits, but only $43\%$ for
10000 digits.

The running times of the serial implementation of \cm, however, reflect
closely the quartic complexity; the computations for the 10000 digit
number are expected to take close to a week and were not carried out.
When compared to the parallel version, it also becomes apparent how the
latter one profits from part of the additional computing power not for
decreasing the wall-clock time, but for shortening the certificates.

The FastECPP implementation in \cm\ has been adopted by the primality
proving community, which has used it for all ECPP records since the first
release of the code in \cm~0.4.0 in May 2022, that is, for all but two
out of the 20 entries at \\
\centerline {
\url {https://t5k.org/top20/page.php?id=27}.
}

Of the AKS/cyclotomy type competitors, the one with the best complexity
known to date is the (probabilistic) algorithm in \cite {Bernstein07}.
It computes an $(N d)$-th power in the ring
\[
\big( (\Z / N \Z)[Y] / (f (Y)) \big) [X] / (X^e - r (Y))
\]
with $f$ of degree
$d \in (\log L)^{O (\log \log \log L)} \subseteq L^{o (1)}$ and
$e \in L^{2 + o (1)}$,
in time $L^{4 + o (1)}$.
While its sequential complexity is comparable to that of FastECPP,
it does not seem possible to reach a wall-clock time complexity of
$L^{3 + o (1)}$ by parallelising it to $\Omega (L)$ cores,
so that FastECPP appears to remain the algorithm of choice for
proofs of large general primes in a parallel setting.

\begin {sloppypar}
\printbibliography
\end {sloppypar}

\end {document}